\documentclass[a4paper,landscape]{article}

\usepackage{a4} 
\usepackage{amsfonts}
\usepackage{amssymb}
\usepackage{amsmath}
\usepackage{amscd}
\usepackage{amsthm}
\usepackage{epsfig}
\usepackage[all]{xy}
\usepackage{color}
\usepackage{bbm}
\usepackage{graphicx}
\usepackage{psfrag}
\usepackage{hyperref}

\newcommand{\R}{\mathbb{R}}

\newcommand{\Z}{\mathbb{Z}}

\newcommand{\Fe}{{\cal F}}
\newcommand{\Ge}{{\cal G}}

\newcommand{\Tot}{\rm Tot}

\newcommand{\ra}{\rightarrow}

\newcommand{\bprf}{\begin{proof}[Proof]}
\newcommand{\eprf}{\end{proof}}

\newtheorem{thm}{Theorem}

\newtheorem{prop}[thm]{Proposition}
\newtheorem{lem}[thm]{Lemma}
\newtheorem{dfn}[thm]{Definition}

\begin{document}
\title{On a topological fractional Helly theorem}
\author{{\Large Stephan Hell}
\thanks{This research was supported by the
    Deutsche Forschungsgemeinschaft within the European graduate
    program `Combinatorics, Geometry, and Computation' (No. GRK
    588/2).}} 
\date{Institut f\"ur Mathematik, MA 6--2,
TU Berlin,\\ D--10623 Berlin, Germany,
hell@math.tu-berlin.de}

\maketitle

%%%%%%%%%%%%%%% abstract %%%%%%%%%%%%%%%%%%
\begin{abstract}We prove a new fractional Helly theorem for families
  of sets obeying topological conditions.  More precisely, we show
  that the nerve of a finite family of open sets (and of subcomplexes
  of cell complexes) in $\R^d$ is $k$-Leray where $k$ depends on the
  dimension $d$ and the homological intersection complexity of the
  family. This implies fractional Helly number $k+1$ for families
  $\Fe$: For every $\alpha>0$ there is a $\beta(\alpha)>0$ such that
  for sets $F_1,F_2,\ldots,F_n\in\Fe$ with $\bigcap_{i\in I}F_i$ for
  at least $\lfloor\alpha{n\choose k+1}\rfloor$ sets
  $I\subseteq\{1,2,\dots ,n\}$ of size $k+1$, there exists a point
  which is common to at least $\lfloor\beta n\rfloor$ of the $F_i$.
  Moreover, we obtain a topological $(p,q)$-theorem. Our result
  contains the $(p,q)$-theorem for good covers of Alon, Kalai,
  Matou\v{s}ek, and Meshulam
  \cite{alon03:_trans_number_hyper_arisin_geomet} as a special case.
  The proof uses a spectral sequence argument. The same method is then
  used to reprove a homological version of a nerve theorem of Bj\"orner.
\end{abstract}
%%%%%%%%%%%%%%% Intro %%%%%%%%%%%%%%%%
\begin{section}{Introduction}\label{sec-intro}
  Helly's theorem is a classical theorem in convex geometry:
  For every finite family $\Fe$ of convex sets in $\R^d$ in which
  every $d$ or fewer sets have a common point we have $\bigcap \Fe
  \not=\emptyset$. Numerous Helly-type results are known; see
  \cite{eckhoff93:_helly_radon_carat} for a survey.
  
  In this paper we are mainly concerned with topological conditions
  for fractional Helly theorems.  Recently new fractional Helly
  theorems have been derived using different approaches; see
  \cite{alon03:_trans_number_hyper_arisin_geomet}, \cite{barany03:_helly},
  \cite{matousek04:_bound_vc_helly}. A finite
  or infinite family $\Fe$ of sets has {\it fractional Helly number
    $k$} if for each $\alpha\in (0,1]$ there is a $\beta(\alpha )>0$
  such that following implications holds: For all $F_1,F_2,\ldots,
  F_n\in\Fe$ such that $\bigcap_{i\in I} F_i\not=\emptyset$ for at
  least $\lfloor\alpha{n\choose k}\rfloor$ index sets $I\in{[n]\choose
    k}$, there exists a point which is in at least $\lfloor\beta
  n\rfloor$ of the sets $F_i$. Here $[n]$ is short for the set
  $\{1,2,\ldots,n\}$, and ${X\choose k}$ for the set of $k$-element
  subsets of a set $X$. 
  There is two main tasks concerning fractional
  Helly theorems: 
  \begin{itemize}
  \item For a family $\Fe$ find $\beta(\alpha)>0$
  optimal, as large as possible. 
  \item Determine new families of sets
  that admit a fractional Helly theorem. What is their fractional
  Helly number?
  \end{itemize}
  In this paper, we focus on the second problem
  motivated by a question of Kalai and Matou\v{s}ek: {\it Is there a
    homological analog of VC-dimension?} We give a positive answer
  to this question in Theorem \ref{thm-fh-acyclic} where homological
  conditions imply a fractional Helly theorem.
  
  The original fractional Helly theorem for convex sets in $\R^d$ by
  Katchalski and Liu can then be stated in the following way.
\begin{thm}[Fractional Helly theorem for convex sets
  \cite{katchalski79}]\label{thm-fh-convex} For every $d\geq 1$, the
  family of convex sets in $\R^d$ has fractional Helly number $d+1$.
\end{thm}
This result was generalized by Alon et
al.~\cite{alon03:_trans_number_hyper_arisin_geomet} to {\it good
  covers} in $\R^d$ where a finite family of sets in whose members are
either all open or all closed, is called a {\it good cover} if
$\bigcap\Ge$ is either empty or contractible for all subfamilies
$\Ge\subseteq\Fe$. Moreover, it was shown by Kalai
\cite{kalai84:_inter} that the optimal  
$\beta(\alpha)$ equals $1-(1-\alpha)^{1/(d+1)}$, so the case $\alpha=1$ implies
Helly number $d+1$. Their proof uses the nerve theorem and the
following proposition for $d^*$-Leray families.
\begin{prop}[Fractional Helly theorem for Leray families 
  \cite{alon03:_trans_number_hyper_arisin_geomet}]\label{prop-fh-leray}
  Let $\Fe$ be a finite $d^*$-Leray family, and let $\Fe^\cap$ the
  family of all intersections of the sets of $\Fe$. Then $\Fe^\cap$
  has fractional Helly number $ d+1$. Moreover, one can take
  $\beta(\alpha)=1-(1-\alpha)^{1/(d+1)}$.
\end{prop}
Here a family $\Fe$ is called {\it $d^*$-Leray} if its nerve complex
$N(\Fe)$ is $d$-Leray: For all induced subcomplexes $L\subseteq
N(\Fe)$ the homology groups $H_n(N(\Fe))$ vanish for all $n\geq d$.
As a special case the fractional Helly property holds for the family
$\Fe$. Proposition \ref{prop-fh-leray} plays a key role in proving one
of the main results of
\cite{alon03:_trans_number_hyper_arisin_geomet}, a $(p,q)$-theorem
for finite good covers.

In this paper we identify topological/homological conditions for
families of sets which imply that their nerve complex is $k$-Leray,
where $k$ depends on the above conditions.  Having Proposition
\ref{prop-fh-leray} in mind we extend the fractional Helly theorem for
good covers to families with higher topological intersection
complexity; see Figure \ref{fig-structure} for the relations between
the results.
\begin{thm}[Topological fractional Helly
  theorem]\label{thm-fh-acyclic} Let $\Fe$ be a finite family of open
  sets (or of subcomplexes of a cell complex) in $\R^d$, and $k\geq d$
  such that for all subfamilies $\Ge\subseteq\Fe$ one of the following
  conditions holds:
\begin{enumerate}\item $\bigcap\cal G$ is empty, or
\item the reduced homology groups of $\bigcap\cal G$ vanish in
  dimension at least $k-|G|$, that is
\[\tilde{H}_n({\bigcap\cal G})=0\text{ for all }n\geq k-|\Ge |.\]
\end{enumerate}
Then $\Fe^\cap$ has fractional Helly number $ k+1$. Moreover, we can
choose $\beta(\alpha)=1-(1-\alpha)^{1/(k+1)}$.
\end{thm}
As in the case of good covers this implies fractional Helly number
$k+1$ for the family $\Fe$.  We call a family of sets as in Theorem
\ref{thm-fh-acyclic} satisfying conditions (i) and (ii) a {\it
  $(k-|\Ge|)$-acyclic} family. 
\begin{figure}[h!]  
\begin{tabular}{l|l|l}\label{tab-cond}
 & good cover & $(k-|\Ge|)$-acyclic \\ \hline&& \\  
$\tilde{H}_n(F_i)$ & $0$\, for all $n\geq 0$ & 
$\left\{\begin{array}{ll} 
0 & \text{ for all } n\geq k-1,\\
\text{arbit.} & \text{ for  $n=0,1,\ldots,k-2$}.
\end{array} 
\right.$\\[3ex]
$\tilde{H}_n(F_{i_1}\cap F_{i_2})$ & $0$\, for all $n\geq 0$ & $\left\{\begin{array}{ll} 
0 & \text{ for all } n\geq k-2,\\
\text{arbit.} & \text{ for  $n=0,1,\ldots,k-3$}.
\end{array} 
\right.$\\
\hspace*{2ex}$\vdots$ & \hspace*{2ex}$\vdots$ &
\hspace*{2.8ex}$\vdots$\\[1ex]
$\tilde{H}_n(F_{i_1}\cap F_{i_{k-1}})$ & $0$\, for all $n\geq 0$ & $\left\{\begin{array}{ll} 
0 & \text{ for all } n\geq 1,\\
\text{arbit.} & \text{ for  $n=0$}.
\end{array} 
\right.$\\[3ex]
$\tilde{H}_n(F_{i_1}\cap\cdots\cap F_{i_t})$ & $0$\, for all $n\geq 0$
&\hspace*{3.1ex}$0$\hspace*{7.2ex}for all $n\geq 0$ and $t\geq k$.
\end{tabular}
\caption{\label{tab-good-acyclic}Homological conditions for good
  covers and for $(k-|\Ge|)$-acyclic families}
\end{figure}

For $k=d$ we obtain fractional Helly
number $k+1$ in a more general setting than good covers. Table
\ref{tab-good-acyclic} shows in the first column
the conditions in the good cover case: All non-empty intersections are
contractible, so their homology vanishes in all dimensions. In the second
column the conditions in the $(k-|\Ge|)$-acyclic case are shown:
Non-empty intersections of $i<k$ sets can have arbitrary homology
groups in dimension less or equal than $k-i-1$.

The case $k>d$ admits even more general families of sets in $\R^d$.
The price one has to pay for increasing the intersection complexity of
the $F_i$ is a higher fractional Helly number. See Figure
\ref{fig-acycli} for an example of the intersection pattern of $k=3$
sets of a $(3-|\Ge|)$-acyclic family in $\R^2$. There the $F_{i_j}$ can
have an arbitrary number of $0$- and of $1$-dimensional holes. The
intersection of two elements $F_{i_j}\cap F_{i_k}$ of our family still
can have an arbitrary number of 0-dimensional holes.

Matou\v{s}ek showed a fractional Helly theorem for families with
bounded VC-dimension; see \cite{matousek04:_bound_vc_helly}.
Matou\v{s}ek's results is not a special case of our result.  Bounded
VC-dimension does not guarantee any Helly property, e.~g.~the family
$\{[n]\setminus\{i\}\,|\, i\in [n]\}$ has bounded VC-dimension, but
no Helly property. Another important example of families with bounded
VC-dimension is the family of all semialgebraic subsets in $\R^d$ of
bounded description complexity. Let's look at a concrete example:
Define a semialgebraic set $F_i=\{x\in\R^d\,|\,x_1^2+x_2^2-i\geq 0\}$,
then the family $\{F_1,F_2,\ldots,F_n\}$ has fractional Helly number
$d+1$. However, we have $\tilde{H}_1(\bigcap_I F_i)=\Z\not=0$ for all
index sets $\emptyset\not=I\subseteq [n]$.

B\'{a}r\'{a}ny and Matou\v{s}ek showed in \cite{barany03:_helly}
that the family of convex lattice sets has fractional Helly number
$d+1$, using a Ramsey-type argument. We can not hope to obtain the
same fractional Helly number using our results as its Helly number is
known to be $2^d$.
\begin{figure}[bht]
\begin{center}
\psfrag{b}{$F_{i_1}$}
\psfrag{a}{$F_{i_2}$}
\psfrag{c}{$F_{i_3}$}
\includegraphics{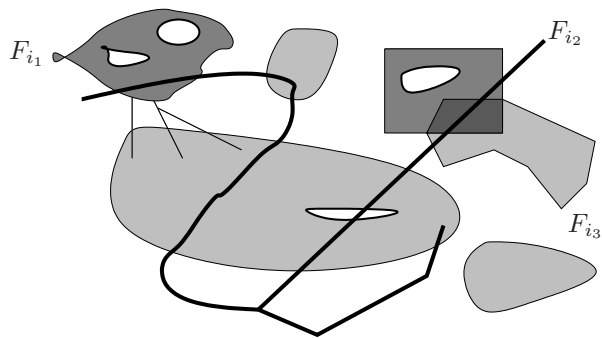}
\end{center}
\caption{\label{fig-acycli}Example of $k=3$
sets of a $(3-|\Ge|)$-acyclic family in $\R^2$}
\end{figure}

The $(p,q)$-theorem for convex sets was conjectured by Hadwiger and
Debrunner, and proved by Alon and Kleitman
\cite{alon92:_pierc_hadwig_debrun}.  For this let $p,q,d$ be integers
with $p\geq q\geq d+1\geq 2$. Then there exists a number HD$(p,q,d)$
such that the following holds: Let $\Fe$ be a finite family of convex
sets in $\R^d$ satisfying the $(p,q)$-condition; that is, among any
$p$ sets of $\Fe$, there are $q$ sets with a non-empty intersection.
Then $\tau(\Fe)\leq \text{HD}(p,q,d)$, where $\tau(\Fe)$ denotes the
{\it transversal number} of $\Fe$, i.~e.~the smallest cardinality of a
set $X\subseteq \bigcup\Fe$ such that $F\cap X\not=\emptyset$ for all
$F\in\Fe$.  It was observed in
\cite{alon03:_trans_number_hyper_arisin_geomet} that the crucial
ingredient in the proof is a fractional Helly theorem for
$\Fe^\cap$. Therefore Theorem \ref{thm-fh-acyclic} implies
immediately a new
$(p,q)$-theorem using the general tools developed in
\cite{alon03:_trans_number_hyper_arisin_geomet}.
\pagebreak
\begin{thm}[$(p,q)$-theorem for $(k-|\Ge|)$-acyclic
  families]\label{thm-pq-acyclic} The assertions of the
  $(p,q)$-theorem also hold for finite $(k-|\Ge|)$-acyclic families of
  open sets (or of subcomplexes of a cell complex) in $\R^d$ where
  $p\geq q\geq k\geq d+1\geq 2$.
\end{thm}
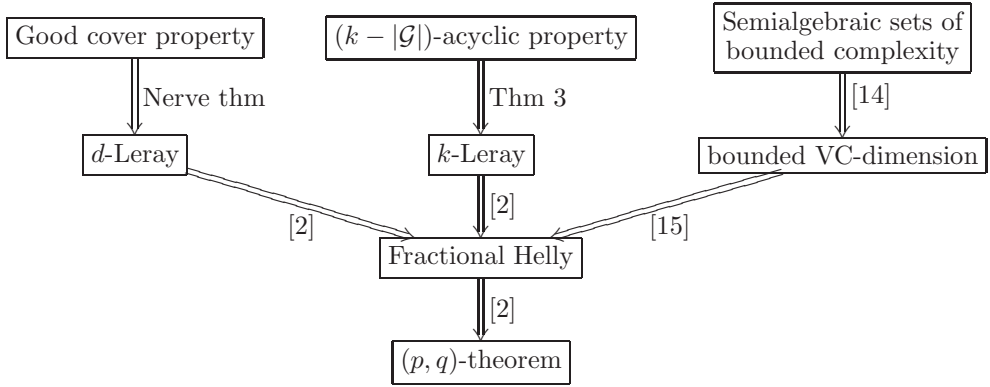
\begin{figure}[ht]
\xymatrix{
*+[F]{\txt{Good cover property}}\ar@{=>}[d]^{\txt{Nerve thm}}
&*+[F]{\txt{$(k-|\Ge|)$-acyclic property}}\ar@{=>}[d]^{\txt{Thm \ref{thm-fh-acyclic}}}
&*+[F]{\txt{Semialgebraic sets of \\bounded complexity}}\ar@{=>}[d]^{\txt{\cite{matousek02:_lectur_discr_geomet}}} \\
*+[F]{\txt{$d$-Leray}}\ar@{=>}[dr]_{\txt{\cite{alon03:_trans_number_hyper_arisin_geomet}}}
&*+[F]{\txt{$k$-Leray}} \ar@{=>}[d]^{\txt{\cite{alon03:_trans_number_hyper_arisin_geomet}}}
&*+[F]{\txt{bounded VC-dimension}} \ar@{=>}[dl]^{\txt{\cite{matousek04:_bound_vc_helly}}}\\
&*+[F]{\txt{Fractional Helly}} \ar@{=>}[d]^{\txt{\cite{alon03:_trans_number_hyper_arisin_geomet}}}\\
&*+[F]{\txt{$(p,q)$-theorem}}
}
\caption{\label{fig-structure}Diagram of the abstract machinery}
\end{figure}
The proof of Theorem \ref{thm-fh-acyclic} uses a spectral sequence
argument. Section \ref{sec-prel} comes with a crash course on
spectral sequences. In Section \ref{sec-acyclic} we prove Theorem
\ref{thm-fh-acyclic}. The same technique is used in Section
\ref{sec-nerve} to reprove a homological versions of a nerve
theorem of Bj\"orner.
\end{section}
%%%%%%%%%%%%%% Preliminaries %%%%%%%%%%%%%%%%%%%%%%%%%%%%%
\begin{section}{Preliminaries}\label{sec-prel}
  As a start we give some background and fix our notation. Spanier's
  book \cite{h.spanier66:_algeb_topol} is a good reference for
  (algebraic) topology.  A topological space $X$ is {\it connected} if
  it is not the disjoint union of two non-empty open subsets. If $X$
  is a connected space then one has $H_0(X,G)=G$ for singular homology
  with coefficients in $G$. A topological space $X$ is ${\it
    contractible}$ if the identity $i:X\ra X$ is homotopic to a constant
  map $c:X\ra X$. The singular homology of a contractible space vanishes in all
  dimensions except for dimension 0 where it equals the coefficient group
  $G$.
  
  Let $\Delta^p$ be the standard $p$-dimensional simplex with vertex
  set $[p+1]$.  The {\it nerve $N(\Fe)$} of a family of sets $\Fe$ is
  the simplicial complex with vertex set $\Fe$ whose simplices are all
  $\sigma\subseteq\Fe$ such that $\bigcap_{F\in\sigma}
  F\not=\emptyset$. For a family $\Fe=\{F_i\,|\, i\in I\}$ of
  subspaces we define the group of singular $n$-chains:
  \[S_n\{\Fe\}=\Z^{\{\sigma:\Delta^p\ra\bigcup\Fe\,|\,\text{im}(\sigma)\subseteq F_i \text{ for some } i\}}\]
  For finite families of open sets (or of subcomplexes of a cell
  complex) the inclusion of the singular chain groups
  $S_*\{\Fe\}\hookrightarrow S_*(\bigcup\Fe)$ induces an isomorphism
  in homology. In the following we write $H_*(X):=H_*(X,\Z)$ for the
  singular homology with integer coefficients of a topological space
  $X$. For simplicity we use for non-empty spaces $X$ also the reduced
  singular homology groups $\tilde{H}_n(X)$, thus saving extra
  considerations for the case $n=0$.  Most of our work also holds for
  arbitrary (co)homology theories and arbitrary coefficients.  As this
  paper is mainly addressed to people who use algebraic topology as a
  tool box we abstain from a more general formulation.
  
  Singular homology however shows anomalies first noticed in
  \cite{barrat62}. Subspaces $A\subseteq\R^d$ which are {\it not nice}
  can have non-vanishing homology in infinitely many dimensions.  For
  this let $A$ be the union of countable many spheres of fixed
  dimension $r>1$ all having one point in common with their diameter
  going to zero, then $A$ is such subspace which is not nice. In the
  literature this example with $r=1$ is also known as the {\it
    Hawaiian earring}. To exclude such not-nice phenomena we consider
  only families $\Fe$ of open sets in $\R^d$, and of subcomplexes of 
  CW-complexes (cell complexes) in $\R^d$. In both cases one has
  $H_n(\bigcup\Fe)=0$ for all $n\geq d$.
  
  {\bf Spectral sequence of a double complex.} Spectral sequences are
  not standard tools in combinatorics so we repeat some definitions
  from \cite{mccleary01:_guide_to_spect_sequen} on {\it spectral sequences
  for homology}; see also
  \cite{basu03:_differ_bound_differ_betti_number_semi} for a short
  introduction. Let $C_{*,*}$ a double complex with two differentials
  $\partial^I:C_{p,q}\ra C_{p-1,q}$ and
  $\partial^{II}:C_{p,q}\ra C_{p,q-1}$ such that
  $\partial^I\partial^{II}+\partial^{II}\partial^I=0$. We associate to
  a double complex its total complex $\Tot (C)_n=\bigoplus_{p+q=n}C_{p,q}$
  with differential $d:=\partial^I+\partial^{II}$. The above relation
  on $\partial^I$ and $\partial^{II}$ implies $d\circ d=0$. In this
  paper we are interested in first quadrant sequences so $C_{p,q}=0$
  for $p<0$ or $q<0$.  Before going into more details we state the
  following tool for the homology $H_*(\Tot(C),d)$.
\begin{thm}[Spectral sequence of a double complex; {\cite[Theorem
  2.15]{mccleary01:_guide_to_spect_sequen}}]\label{thm-double-complex}
  Given a double complex $(C_{*,*},\partial,\tilde{\partial})$ there
  are two spectral sequences $(E^r_{*,*}, d^r)$ and $(
  \tilde{E}^r_{*,*}, \tilde{d}^r)$ with \[ E^2_{*,*}\cong
  H^{\partial}_{*,*}H^{\tilde{\partial}}(C)\;\;\text{ and }\;\;
  \tilde{E}^2_{*,*}\cong H^{\tilde{\partial}}_{*,*}H^{\partial}(C). \]
  If $C_{p,q}=0$ for $p<0$ and $q<0$ then both spectral sequences
  converge to $H_*(\Tot (C),d)$.
\end{thm}
Here $H^\partial_*(C)$ stands for the homology of $C_{*,*}$ with respect to
the boundary $\partial$. The boundary $\tilde{\partial}$ induces a
boundary on $H^\partial_*(C)$ so that
$H^{\tilde{\partial}}_{*,*}H^{\partial}(C)$ is well defined. We now
repeat in detail the construction of the
spectral sequences to a given double complex.
\begin{dfn}[Spectral sequence of homological type] \label{def-spectral-sequence}
  A spectral sequence is a collection of differential
  bigraded modules, that is, for $r=1,2,3,\ldots$, and $p,q\geq 0$ we
  have a module $E^r_{p,q}$, and furthermore differentials
  $d^r:E^r_{p,q}\ra E^r_{p-r,q+r-1}$. Finally,
  $E^{r+1}_{*,*}\cong H(E^r_{*,*},d_r)$ for all $r\geq 1$.
\end{dfn}
\noindent Figure \ref{fig-differentials} shows the differentials
$d^r$; for $r$ big enough $d^r$ hits the row 0.
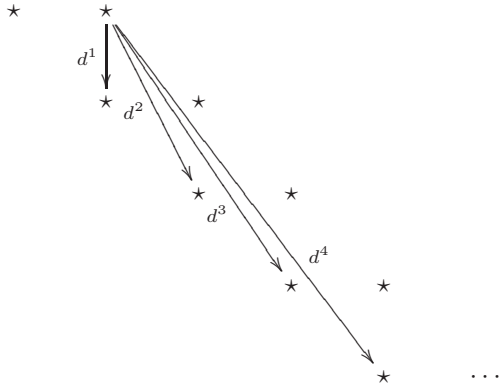
\begin{figure}[hbt]
\xymatrix{
   %\star & & & & & &\\
   \star & \star \ar[d]_{d^1} \ar[ddr]_{d^2}\ar[dddrr]_(.7){d^{3}}
   \ar[ddddrrr]^(.7){d^4}& & & & & \\
   &\star &\star & & & & \\
   & & \star &\star & & &\\
   %& & & \ddots & \ddots & & \\
   & & & \star & \star  & \\
   & & & & \star &  \cdots\\
 }
\caption{\label{fig-differentials}The differentials $d^r:E^r_{p,q}\ra E^r_{p-r,q+r-1}$.}
\end{figure}
Let $(C_{*,*},\partial,\tilde{\partial})$ be a double complex with
$\partial:C_{p,q}\ra C_{p-1,q}$, $\tilde{\partial}:C_{p,q}\ra
C_{p,q-1}$, and total differential $d=\partial + \tilde{\partial}$. We
define two filtrations of the total complex $(\Tot (C)_n,d)$ as follows: \[
F_m(\Tot (C)_n)=\bigoplus_{p\leq m}C_{p,n-p}\,\,\text{ and }\,\,
\tilde{F}_m(\Tot (C_n))=\bigoplus_{q\leq m}C_{n-q,q}. \] From now on
we do
construction for the first filtration. This goes over to
the construction for the second filtration by reindexing the double
complex as its transpose $C^T_{p,q}=C_{q,p},\,
\partial^T=\tilde{\partial},\, \tilde{\partial}^T=\partial$.  The
filtration is increasing and respects the total differential $d$:
\[ 0\cdots\subseteq F_{m-1}(\Tot (C)_n) \subseteq F_m(\Tot (C)_n)\subseteq
\cdots \Tot (C)_n\]
and,
\[d(F_m(\Tot(C)_n))\subseteq F_m(\Tot(C)_n).\]
Since $d$ respects the filtration the homology of $(C,d)$ inherits a
filtration
\[0\subseteq\cdots\subseteq F_{m-1}( H(\Tot(C),d))\subseteq F_m(H(\Tot(C),d))
\subseteq\cdots\subseteq H(\Tot(C),d),\] where
$F_m(H(\Tot(C),d)):=\iota^m_*(H(F_m(\Tot(C)),d))\subseteq
H(\Tot(C),d)$, and \linebreak $\iota^m:F^m(\Tot(C))\hookrightarrow
\Tot(C)$ is the inclusion.  We know from Theorem \ref{thm-double-complex}
that the associated spectral sequence converges to $H(\Tot(C),d)$;
more precisely, there is a $r>0$ such that $d_r=0$ is trivial. This is
called {\it the sequence collapses at term $r$}.  Therefore we have
\[E^{r+1}_{p,q}=E^\infty_{p,q}\cong
F^p(H_{p+q}(\Tot(C),d))/F^{p-1}(H_{p+q}(\Tot(C),d)). \] Hence, one
gets $H_n(C,d)\cong\bigoplus_{p+q=n} E_{p,q}^\infty$ in the case of
homology with field coefficients. In the case of integer coefficients
this leads to an extension problem. In this paper the extension
problem is trivial as only one of the groups $E_{p,q}^\infty$ with
$n=p+q$ is different from zero, namely either $E_{n,0}^\infty$ or
$E^\infty_{0,n}$. Therefore $H_n(\Tot(C),d)$ equals either
$E_{n,0}^\infty$ or $E^\infty_{0,n}$.

Consider the following definitions
for $r\geq 0$
\begin{eqnarray*}
Z^r_{p,q} & = & F_p(\Tot(C)_{p+q})\cap d^{-1}(F_{p-r}(\Tot(C)_{p+q-1}))  \\
B^r_{p,q} & = & F_p(\Tot(C)_{p+q})\cap d(F_{p+r}(\Tot(C)_{p+q+1})) \\
Z^\infty_{p,q} & = & F_p(\Tot(C)_{p+q}) \cap \text{ker}(d) \\
B^\infty_{p,q} & = & F_p(\Tot(C)_{p+q}) \cap \text{im}(d) 
\end{eqnarray*}
Then this leads to a tower of submodules
\[ B^0_{p,q}\subseteq B^1_{p,q}\subseteq\cdots\subseteq B^\infty_{p,q}
\subseteq Z^\infty_{p,q}\subseteq\cdots\subseteq Z^1_{p,q}\subseteq
Z^0_{p,q}.\] For a first quadrant sequence we get that for
$r>\max\{p,q\}$\[ Z^r_{p,q} = Z^\infty_{p,q}\,\,\,\text{, and }\,\,\,
B^r_{p,q}=B^\infty_{p,q}\] hold. This insures the convergence of our
spectral sequence. Define \[ E^r_{p,q} = Z^r_{p,q}/ (Z^{r-1}_{p-1,q}+
B^{r-1}_{p,q}) \] It can be checked that the differential
$d:Z^r_{p,q}\ra Z^r_{p-r,q+r-1}$ induces a unique homomorphism
$d^r$ such that the following diagram commutes: \[\xymatrix{
  Z^r_{p,q}\ar[r]^d\ar[d]_{\eta^{r}_{p,q}} & Z^r_{p-r,q+r-1}
  \ar[d]^{\eta^{r}_{p-r,q+r-1}} \\
  E^r_{p,q}\ar[r]^{d^r} & E^r_{p-r,q+r-1} }\] where the maps
$\eta^r_{p,q}:Z^r_{p,q}\ra E^r_{p,q}$ are the canonical
projections. The same constructions leads to the spectral sequence
$(\tilde{E}^r,\tilde{d}^r)$, $\tilde{d}^r:\tilde{E}^r_{p,q}\ra
\tilde{E}^r_{p+r-1,q-r}$, for the second filtration.
\end{section} 

%%%%%%%%%%%% Proof of the main result
\begin{section}{On $(k-|\Ge|)$-acyclic families}\label{sec-acyclic}
  In this section we prove our main result Theorem
  \ref{thm-fh-acyclic}. The following lemma is the key argument in
  extending the fractional Helly theorem from good covers to
  $(k-|\Ge|)$-acyclic families.
\begin{lem}\label{lem-acyclic-homology}For $k\geq 0$ let $\Fe$ be a finite
  $(k-|\Ge|)$-acyclic family of open sets (or of subcomplexes of a cell
  complex) in $\R^d$. 
Then $H_n(\bigcup\Fe)\cong H_n(N(\Fe))$ for all $n\geq k$.
\end{lem}
To prove Lemma \ref{lem-acyclic-homology} we first define a suitable
double complex $(C_{*,*},\partial,\tilde{\partial})$.  Then we compute
its $E^2_{*,*}$- and $\tilde{E}^2_{*,*}$-term which are shown in
Figures \ref{fig-e2-term} and \ref{fig-ee2-term}. 
Finally, we apply Theorem \ref{thm-double-complex}
to get the conclusion.
\bprf 
For $\Fe=\{F_1,F_2,\ldots,F_n\}$ define a double complex
\[ C_{p,q} := \bigoplus_{J\subseteq [n],\,|J|=q+1} S_p(\bigcap_{j\in
  J}F_j)\,\,\text{ for }p,q\geq 0.\] Let $\partial:=\oplus\,
\partial:C_{p,q}\ra C_{p-1,q}$ be the usual singular boundary
operator. For index sets $J' \subseteq J\subseteq [n]$ let
$i^{J,J'}:\bigcap_{j\in J}F_j\ra\bigcap_{j\in J'}F_j$ be the
inclusion. The group $S_p(\bigcap_{j\in J}F_j)$ is freely generated by
the set of singular $p$-simplices $\sigma:\Delta^p
\ra\bigcap_{j\in J}F_j$. Define $\tilde{\partial}$ component-wise 
on the elements $c=\sum r_\sigma \sigma$ of $S_p(\bigcap_{j\in
  J}F_j)$:\[ c \mapsto \tilde{\partial}(\sigma) :=
(-1)^p\sum_\sigma\sum_{i=0}^q(-1)^i r_\sigma i^{J,J_i}_*(\sigma)
\,\,\in\bigoplus_{J\subseteq [n],|J|=(q-1)+1} S_p(\bigcap_{j\in J}F_j),\]
where $J_i=\{j_0<j_1<\ldots<\hat{j_i}<\dots <j_q\}$ is the set
obtained from $J$ by deleting the element $j_i$, and the factor
$(-1)^{p}$ is added to guarantee $\partial\tilde{\partial} +
\tilde{\partial}\partial=0$.

We show that (i)
\[E^2_{p,q}=H^{\partial}_{p,q}H^{\tilde{\partial}}(C)=
\left\{ \begin{array}{ll}
    H_p(\bigcup \Fe) & \text{for } q=0, \\
    0 & \text{otherwise,}
\end{array}
\right. \] and that (ii)
\[\tilde{E}^2_{p,q}=H^{\tilde{\partial}}_{p,q}H^{\partial}(C)= 
\left\{ \begin{array}{ll} 
H_q(N(\Fe)) & \begin{array}{l}\text{for }p=0,\, q\geq k\end{array} \\
0 & \begin{array}{l}\text{for }p=i,\,q=n-i
\text{, for all }\\ 1\leq i\leq n\text{ and }n\geq k-1.
\end{array}
\end{array}
\right.\] 
\begin{figure}[htb]
\xymatrix@!C{
  H_n(\bigcup\Fe) & 0   & 0   & 0  & \cdots\\
  \vdots & 0 & 0 & 0 & \cdots\\
  H_1(\bigcup\Fe) & 0 & 0  & 0  & \cdots \\
  H_0(\bigcup\Fe) & 0 & 0 & 0 & \cdots 
}
\caption{\label{fig-e2-term}$E^2$-term in the proof in the proof of Lemma \ref{lem-acyclic-homology}}
\end{figure}
The $E^2$-term is equal to zero except for the zeroth column where we
get the homology of $\bigcup \Fe$, see Figure \ref{fig-e2-term}. The
$\tilde{E}^2$-term contains in the zeroth row the homology of $N(\Fe)$
for dimensions greater or equal than $k$, at the same time vanish all
terms above the $(k-1)$-th anti-diagonal, see also Figure
\ref{fig-ee2-term}. Hence the $\tilde{E}^r$-sequence collapses at term
2 above this anti-diagonal. Using Theorem \ref{thm-double-complex} we
obtain
\[ H_n(\bigcup\Fe)\cong H_n(\Tot(C))\cong H_n(N(\Fe))\text{ for all }
n\geq k.\]

To obtain (i) we first compute $H^{\tilde{\partial}}(C_{p,*})$.  Using
the definition of singular homology one gets:
\[ C_{p,q}=\bigoplus_{J\subseteq [n],|J|=q+1}
\;\,\bigoplus_{\sigma:\Delta^p\ra
  X,\text{im}(\sigma)\subseteq\bigcap_{j\in J} F_j}\Z\] For
$\sigma:\Delta^p\ra X$ let $J_\sigma$ be the maximal subset
$J\subseteq [n]$ with $\text{im}(\sigma)\subseteq \bigcap_{j\in J}
F_j$. Then this leads to:
\[ C_{p,q}= \bigoplus_{\sigma:\Delta^p\ra X}\;\,
\bigoplus_{J\subseteq [n],\,J\subseteq J_\sigma,\,|J|=q+1} \Z\] Notice
that for $J_\sigma=\emptyset$ a (innocent) zero is added to the direct
sum. The boundary $\tilde{\partial}$ has no effect on $\sigma$ so that
we can look at every component
\[\bigoplus_{J\subseteq [n],\,J\subseteq J_\sigma,\,|J|=q+1} \Z\]
separately. For $J_\sigma\not=\emptyset$ one can check that this
equals the simplicial chain complex of a $(|J_\sigma|-1)$-dimensional
simplex. The simplicial homology of the simplex vanishes in all
dimensions except dimension $0$ where it equals $\Z$. Using that
\[\bigoplus_{\sigma:\Delta^p\ra X}\;\,
\bigoplus_{J\subseteq [n],\,J\subseteq J_\sigma,\,|J|=1}
\Z=S_p(\Fe),\] we obtain:
\[H_q^{\tilde{\partial}}(C_{p,*})=
\left\{ \begin{array}{ll} 
S_p(\Fe) & \text{for }q=0,\\
0        & \text{otherwise.}
\end{array}
\right.    
\]
The induced boundary $\partial$ can be identified with
the usual boundary on $S_*(\Fe)$.
For families of open sets or of subcomplexes of a cell complex 
the inclusion $S_*(\Fe)\ra S_*(X)$ induces an isomorphism
$H_*(\Fe)\ra H_*(X)$ in homology. 

To prove (ii) we start with computing $H^\partial (C_{*,q})$, 
we make use of the $(k-|\Ge|)$-acyclicity of our family. 
By the definition of $\partial$ we know that
\[ H^\partial_p(C_{*,q})={\displaystyle
     \bigoplus_{J\subseteq[n],\,|J|=q+1}H_p(\bigcap_{j\in J}F_j).}
\]
If $\bigcap_{j\in J} F_j=\emptyset$ then the contribution to this sum
equals $0$. In the case $\bigcap_{j\in J} F_j\not=\emptyset$ the
summand depends on our assumption $\tilde{H}_n(\bigcap_{j\in J}F_j)=0$
for $n\geq k-|J|$. E.~g.~for $p=0$ and $|J|\geq k$ the summand equals $\Z$,
and in general this leads to
\[H^\partial_p(C_{*,q})=
\left\{ \begin{array}{cl}{ \displaystyle
      \bigoplus_{J\subseteq[n],\,|J|=q+1,\,\bigcap_{j\in
          J}F_j\not=\emptyset}\Z} 
& \text{for }p = 0,\,q\geq k-1 \\
0 
& \begin{array}{l}\text{for }p=i,\,q=n-i
\text{, for all }\\ 1\leq i\leq n\text{ and }n\geq k-1.
\end{array}
\end{array}
\right. 
\] 
\begin{figure}[htb]
\xymatrix@!C{
  0  & 0   & 0   & 0  & \cdots\\
  \star & 0 & 0      & 0 & \cdots\\
  \star & \star & 0     & 0  & \cdots \\
  \star & \star & \star &  H_{k}(N({\cal G})) & \cdots 
}
\caption{\label{fig-ee2-term}The $\tilde{E}^2$-term in the proof of Lemma \ref{lem-acyclic-homology}}
\end{figure}
It is easy to check that the induced chain complex $(H_0^\partial
(C_{*,*}),\tilde{\partial})$ and the simplicial chain complex
$C_*(N(\Fe))$ are isomorphic in dimension $\geq k$. The chain groups are
even isomorphic in dimensions $\geq k-1$, but the differential differ
in general in dimension $k-1$. To see that differentials coincide up
to a sign in dimensions $\geq k$ remember that $H_0(\bigcap_{j\in J}
F_j)$ is freely generated by the class of any $0$-simplex, and
$i^{J,J_i}_*$ maps the class of a $0$-simplex on a class of a
$0$-simplex. Hence row 0 of the $E^2$-term contains the homology
$H_n(N({\cal G}))$ for $n\geq k$.  
\eprf

After these preparations the proof of Theorem \ref{thm-fh-acyclic} is
short. In the next section we make use of the same double complex
$(C_{*,*},\partial,\tilde{\partial})$. Similar computations lead to a
nice proof of a homological version of a nerve theorem of Bj\"orner.
\bprf (of {\bf Theorem \ref{thm-fh-acyclic}}) We prove that our family
$\Fe$ is $k$-Leray, hence Theorem \ref{thm-fh-acyclic} follows
immediately from Theorem \ref{prop-fh-leray}. For this let $L\subseteq
N(\Fe)$ be an induced subcomplex. Then $L$ is of the form $N(\Ge)$ for
some $\Ge\subseteq\Fe$. The family $\Ge$ is again $(k-|\Ge|)$-acyclic so
that Lemma \ref{lem-acyclic-homology} implies
\[ H_n(\bigcup\Ge)\cong H_n(N(\Fe))\text{ for all }n\geq k.\]
Finally we have $H_n(\bigcup\Ge)=0$ for all $n\geq d$ as
$\bigcup\Ge$ is a {\it nice} subset of $\R^d$. 
\eprf
\end{section}
%%%%%%%%% Nerve Theorems %%%%%%%%%%%%%%%%%%%%%%%%
\begin{section}{On nerve theorems}\label{sec-nerve}
  The nerve theorem is a standard tool in topological combinatorics.
  It was first obtained by Leray \cite{leray45:_sur}; see Bj\"orner
  \cite{bjoerner95:_topol} for a survey on nerve theorems.  Here we
  prove a homological version of a nerve theorem of Bj\"orner.
  For this recall that for $k\geq -1$ a topological space $X$ is {\it
    $k$-connected} if for every $l=-1,0,1,\ldots,k$, each continuous
  map $f:S^l\ra X$ can be extended to a continuous map
  $\bar{f}:B^{l+1}\ra X$. Here the $(-1)$-dimensional sphere $S^{-1}$
  is interpreted as $\emptyset$, and the $0$-dimensional ball $B^0$ as
  a single point. Hence $(-1)$-connected means non-empty, and
  $0$-connected means pathwise connected. 
\begin{thm}[Nerve Theorem II \cite{bjoerner03:_nerves}, homology version]\label{thm-nerve-II}
  For $k\geq 0$ let $\Fe$ be a family of open sets (or a finite family
  of subcomplexes of a cell complex) such that every $\bigcap \Ge$ is
  empty or $(k-|\Ge|+1)$-connected for all non-empty subfamilies
  $\Ge\subseteq\Fe$.  Then $H_n(X)=H_n(N(\Fe))$ holds for
  $X=\bigcup\Fe$ and all $n\leq k$.
\end{thm}
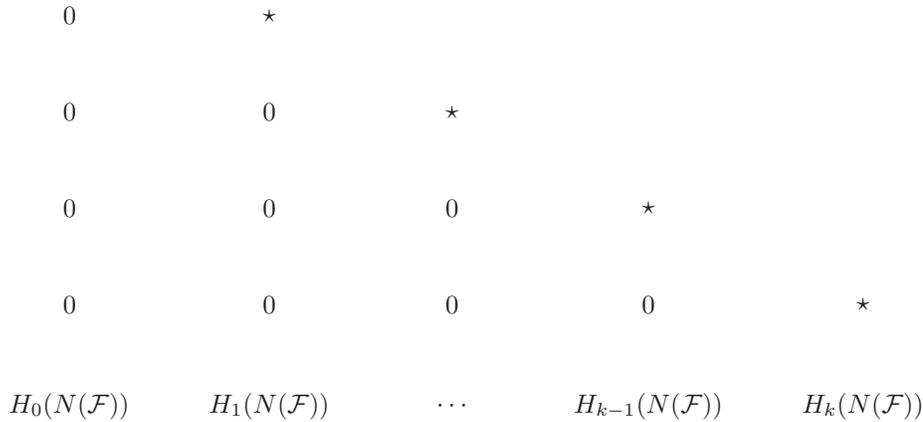
\begin{figure}[thb]
\xymatrix{
  0 &  \star \\
  0 & 0  & \star \\
  0 & 0 & 0 & \star \\
  0 & 0 & 0 & 0 & \star  \\
  H_0(N(\Fe)) & H_1(N(\Fe)) & \hspace*{2.4ex}\cdots\hspace*{2.4ex} & H_{k-1}(N(\Fe)) & H_k(N(\Fe))
}
\caption{\label{fig-nerve}The $\tilde{E}^2$-term in the proof of Theorem \ref{thm-nerve-II}}
\end{figure}
\bprf As in the proof of Lemma \ref{lem-acyclic-homology} we have that
$H_*(X)\cong H_*(\Tot(C),d)$. For a $k$-connected space $X$ we know
from a famous theorem of Hurewicz that $H_n(X)=0$ for all $n\leq k$.
Using the conditions on the connectivity of $\Ge$ and
analogous arguments as in the proof of Lemma
\ref{lem-acyclic-homology} the $\tilde{E}^2$-term looks as in Figure
\ref{fig-nerve}.
Hence we see that $H_n(N(\Fe))\cong H_n(C,d)$ for all $n\leq k$.
\eprf
\end{section}

%%%%%%%%%%%%%%%%%%%%%% Danksagung
\section*{Acknowledgement}
I thank Carsten Schultz for suggesting the use of the spectral
sequence argument, and my advisor G\"unter M.~Ziegler for helpful
discussions. Many thanks also to Ji\v{r}\'{i} Matou\v{s}ek for his
warm hospitality at KAM, Prague, and for introducing me to the
problem.

%%%%%%% Literatur %%%%%%%%%%%%%%%%%%%%%%%%%%%%%%%%%%%%%%%%%%%%%%%

\nocite{chichilnisky93:_inter}
\nocite{alon95:_bound}

\end{document}